# Improved Monte Carlo Variance Reduction for Space and Energy Self-Shielding

S. C. Wilson and R. N. Slaybaugh*

*Bechtel Marine Propulsion Corporation*
*P.O. Box 79, West Mifflin, Pennsylvania 15122-0079*

**Abstract** — *Continued demand for accurate and computationally efficient transport methods to solve optically thick, fixed-source transport problems has inspired research on variance-reduction (VR) techniques for Monte Carlo (MC). Methods that use deterministic results to create VR maps for MC constitute a dominant branch of this research, with Forward Weighted–Consistent Adjoint Driven Importance Sampling (FW-CADIS) being a particularly successful example. However, locations in which energy and spatial self-shielding are combined, such as thin plates embedded in concrete, challenge FW-CADIS. In these cases the deterministic flux cannot appropriately capture transport behavior, and the associated VR parameters result in high variance in and following the plate.*

*This work presents a new method that improves performance in transport calculations that contain regions of combined space and energy self-shielding without significant impact on the solution quality in other parts of the problem. This method is based on FW-CADIS and applies a Resonance Factor correction to the adjoint source. The impact of the Resonance Factor method is investigated in this work through an example problem. It is clear that this new method dramatically improves performance in terms of lowering the maximum 95% confidence interval relative error and reducing the compute time. Based on this work, we recommend that the Resonance Factor method be used when the accuracy of the solution in the presence of combined space and energy self-shielding is important.*

## I. INTRODUCTION

Continued demand for accurate and computationally efficient transport methods to solve optically thick, fixed-source problems has inspired many to research variance-reduction (VR) techniques for Monte Carlo (MC). Methods that use adjoint deterministic results to create VR maps for MC constitute a dominant branch of this research.[1] Among these, the Forward Weighted–Consistent Adjoint Driven Importance Sampling (FW-CADIS) method[2] is particularly useful for fixed-source problems in which the global flux solution is of interest. FW-CADIS can accelerate global fixed-source calculations such that it becomes feasible to use MC to solve them.[1]

The FW-CADIS method requires a forward deterministic calculation to create an adjoint source and an adjoint deterministic calculation to create an importance map for MC. In general, the deterministic solutions do not have to be very accurate to accelerate MC. There may, however, be pathological cases where the deterministic scalar flux cannot appropriately capture behavior, and even the current best methods fail to accelerate the solution in some locations.

In this work we have identified one such pathological case, which occurs when spatial self-shielding and energy self-shielding are present in the same area, for example, neutrons traveling through long, resonance-dominated materials surrounded by moderator. The combination of spatial and energy self-shielding presents impediments to defining correct deterministic multigroup (MG) cross sections and, thus, to effectively reducing MC solution variance. Situations where thin plates containing a resonance material (e.g., $^{56}$Fe) cut through moderator (e.g., water or concrete) are common.

As the path length in the resonance material increases, the spectral shape of the flux begins to mirror the resonance structure. Physically, neutrons scatter into the "dips" of the resonances, and if the scattered direction points directly through the resonance material, their mean free path before interaction may become very large. This results in a spectral peak in the resonance region at the exit of the resonance material. This peak will be substantially underpredicted if the deterministic cross-section generation code uses a smooth flux shape to generate MG cross sections without correcting for this behavior.

---

*E-mail: slaybaugh@berkeley.edu

A popular and often effective method to account for self-shielding in cross sections is the Bondarenko self-shielding correction.[3] Even with this correction, there can be challenges in producing correct deterministic cross-section sets for problems with long, thin resonance materials surrounded by moderator. For FW-CADIS flux optimization, it appears that even correcting the deterministic cross sections with the Bondarenko method may not correct the importance function sufficiently to significantly reduce MC variances in these regions. While a solution with acceptable statistical error is achievable, one must resort to the "brute force" methods of increasing tally volumes or running more histories to acquire acceptable solution statistics.

One theory is that the lack of angle dependence in the FW-CADIS implementation used here, combined with the need to use a coarse MG energy structure in deterministic calculations, results in inadequate VR for resonance energy (hundreds of keV) neutrons that scatter into a flight direction traveling down the plate. In this case, a single particle that would have a high weight without weight control and would increase the statistical uncertainty upon scoring is instead split many times as it travels through the plate. The splitting produces many duplicate particles that all have the same weight. These copies do not scatter or otherwise react out of their original energy-angle state. This results in many (hundreds of thousands in some cases) copies of the same particle being scored on the tally mesh. Because these duplicate particles all have the same value, simulating many copies is unnecessary. If the VR were closer to "ideal," only unique particles would be scored; simulating enough unique particles would reduce the variance.

This work presents results for a representative test problem that illustrates the difficulty of calculating deterministic and MC solutions for geometries with a thin plate embedded in moderator. This work also presents results for the same test problem using a newly developed method that substantially mitigates the associated statistical error problems. This method is based on FW-CADIS but applies a correction factor to the adjoint source. This correction is based on forward fluxes from two different cross-section sets, each generated with different cross-section corrections in the areas with spatial and energy self-shielding. The results show improvement at the exit of the plate without significant impact in the solution quality in the other parts of the problem.

The remainder of this work is organized as follows. Background information used in the new method is provided in Sec. II. Section III illustrates the problem and discusses attempts to resolve the issue with standard methods. The new method is presented in Sec. IV, which is followed by results and discussion in Sec. V and conclusions in Sec. VI.

## II. BACKGROUND

To understand the problem and solution, some background information is helpful. First, FW-CADIS is discussed. This is followed by information about cross-section processing.

### II.A. FW-CADIS

The information in this subsection is largely derived from Ref. 2. The equations are presented as angle independent, though the formulation can be extended to include angle in a straightforward manner.

#### II.A.1. CADIS

Define a response in a volume $V_d$, with a response function $\Sigma_d$, as

$$R = \int_E \int_{V_d} \Sigma_d(\mathbf{r},E)\phi(\mathbf{r},E)dVdE \ , \qquad (1)$$

where $\mathbf{r}$ includes all three spatial dimensions and $\phi(\mathbf{r},E)$ is the scalar neutron flux in neutrons per unit area per unit energy.

In Eqs. (2) through (6), the space- and energy-dependence notations are frequently dropped for brevity. The forward and adjoint transport equations are defined as

$$H\phi = q \quad \text{(forward)} \qquad (2)$$

and

$$H^\dagger \phi^\dagger = q^\dagger \quad \text{(adjoint)} \ , \qquad (3)$$

where

$H$ = transport operator

$q$ = source

† = adjoint quantity.

These obey the identity

$$\langle H\phi, \phi^\dagger \rangle = \langle H^\dagger \phi^\dagger, \phi \rangle \ , \qquad (4)$$

and therefore,

$$\langle q, \phi^\dagger \rangle = \langle q^\dagger, \phi \rangle , \quad (5)$$

where $\langle \cdot \rangle$ is the inner product, which is an integral over volume and energy in this case. If we let $q^\dagger = \Sigma_d$, then

$$\langle q^\dagger, \phi \rangle = \langle \Sigma_d, \phi \rangle = R = \langle q, \phi^\dagger \rangle . \quad (6)$$

Equation (6) expresses that $\phi^\dagger$ represents the expected contribution of a source particle to the response given the source $q$.

If an adjoint transport calculation is performed where the adjoint source $q^\dagger$ is defined as the local response of interest (e.g., $\Sigma_d$), the resulting adjoint flux is effectively a map of the importance of all of the parts of phase space to that response. This physical interpretation can be used to create VR parameters.

Importance values designed to accelerate MC calculations can thus be defined as

$$imp(\mathbf{r},E) = \frac{\phi^\dagger(\mathbf{r},E)}{\langle q(\mathbf{r},E), \phi^\dagger(\mathbf{r},E) \rangle} = \frac{\phi^\dagger(\mathbf{r},E)}{R(\mathbf{r},E)} , \quad (7)$$

where the estimates for $\phi^\dagger$ and $R$ are generated with a coarse deterministic calculation. The weighting of the map by $R$ is intended to ensure consistency with the source. Throughout this work this set of importance values is referred to as the importance map.

### II.A.2. Forward Weighting

If a global quantity is of interest instead of a local quantity, a slightly different approach is in order. In this case, results with uniformly low statistical uncertainty are typically of interest. In 2001, Cooper and Larsen[4] suggested that a uniform particle distribution throughout the phase space would result in a relatively uniform level of statistical uncertainty.

With this in mind, Wagner et al.[2] noted that the Consistent Adjoint Driven Importance Sampling (CADIS) method could be extended to target a nonphysical response of unity throughout the phase space by using a deterministically generated estimate of the forward flux. If, for example, the desired end result is energy- and space-dependent flux everywhere, $q^\dagger$ can be set to the inverse of the forward flux such that

$$q^\dagger(\mathbf{r},E) = \frac{1}{\phi(\mathbf{r},E)} , \quad (8)$$

and then,

$$\begin{aligned} R(\mathbf{r},E) &= \langle q^\dagger(\mathbf{r},E), \phi(\mathbf{r},E) \rangle \\ &= \int_E \int_V \frac{1}{\phi(\mathbf{r},E)} \phi(\mathbf{r},E) dV dE \\ &\approx 1 \quad \forall \mathbf{r}, E . \end{aligned} \quad (9)$$

Targeting other results such as energy-integrated flux or a global response will use different $q^\dagger$s, but the principle is the same.

This extension of CADIS requires a forward deterministic calculation to estimate $\phi$ and an adjoint calculation to estimate $\phi^\dagger$. It is worth noting that practical implementations of FW-CADIS exclude angle because of the difficulty presented by storage and data management associated with using angular flux rather than scalar flux. FW-CADIS is the springboard from which the new method was developed, and the implementation of the new method also excludes angle.

### II.B. Cross-Section Processing

While MC codes sample continuous cross-section data, deterministic codes require the computation of MG cross sections. Despite the fact that the deterministic calculations needed by FW-CADIS do not need to be particularly accurate for the method to be valuable, the accuracy of the cross sections has an impact on the identified pathological case discussed in this work. As a result, the manner in which MG cross sections are processed is an important component of the new method. The following is a brief summary of the Bondarenko method to account for self-shielding and derives from Refs. 3, 5, and 6; see these for details about approximations and assumptions in the Bondarenko method.

In the general case, the microscopic MG cross section for reaction $x$ in nuclide $j$ and group $g$ has units of area and is defined as

$$\sigma_{x,g}^{(j)} = \frac{\langle \sigma_x^{(j)}(u) W(u) \rangle}{\langle W(u) \rangle} , \quad (10)$$

where the inner product is now an integration over the lethargy interval of that group, and $W(u)$ is a lethargy-dependent weighting function. The most basic implementation creates $\sigma_{x,g}^{(j)}$ using $W(u) = \phi_\infty(u)$, which is the infinite medium flux resulting from a fission source in moderator with negligible absorption. This weighting function gives a cross section denoted by $\sigma_{x,g}^{(j)}(\infty)$, termed infinitely dilute. This approach is inaccurate in a variety of cases.

In the Bondarenko method, it is assumed that the differences between the exact form of the spectrum and the assumed weighting function are inversely proportional to the macroscopic total cross section; thus, $W(u) = \phi_\infty(u)(1/\Sigma_t(u))$. When a material contains a mixture of materials that may have resonances, the weighting spectrum is corrected to account for the effect of the other isotopes. The correction requires a background cross section, $\sigma_0$, which is the sum of the total cross sections of all other elements in the medium, each with number density $N$:

$$\sigma_0^{(j)} = \frac{1}{N_j} \sum_{m \neq j} \sigma_t^{(m)} N_m \ . \quad (11)$$

To correct the MG cross sections, change the weighting function to $W(u) = \phi_\infty(u)(1/(\sigma_t^{(j)}(u) + \sigma_0^{(j)}))$, and then define the self-shielding factors as

$$f_{x,g}^{(j)}(\sigma_0^{(j)}) = \frac{1}{\sigma_{x,g}^{(j)}(\infty)} \frac{\left\langle \sigma_x^{(j)}(u) \frac{\phi_\infty(u)}{\sigma_t^{(j)}(u) + \sigma_0^{(j)}} \right\rangle}{\left\langle \frac{\phi_\infty(u)}{\sigma_t^{(j)}(u) + \sigma_0^{(j)}} \right\rangle} \ . \quad (12)$$

These factors are applied to the infinite medium cross sections. Thus,

$$\sigma_{x,g}^{(j)}(\sigma_0^{(j)}) = f_{x,g}^{(j)}(\sigma_0^{(j)}) \times \sigma_{x,g}^{(j)}(\infty) \ . \quad (13)$$

If there is more than one nuclide with resonances in a medium, this can be accounted for by using the corrected total cross sections computed with Eq. (13) in Eq. (11). This discussion excluded temperature dependence for brevity, though the noninfinitely dilute terms are temperature dependent.

When a given nuclide is very dilute in a medium, $\sigma_0^{(j)}$ becomes much larger than $\sigma_t^{(j)}$, and the shape of the weighting function approaches that of the mixture excluding the nuclide; when it is concentrated, $\sigma_0^{(j)}$ becomes small, and the resonances in the concentrated nuclide have a larger impact.

An additional correction is added when a "thin slab" of resonance absorber is placed within a "thick slab" of moderator. The absorbing slab is considered a separate region, and the equations above are modified to account for its interaction with the surrounding moderator. When the resonance slab's thickness is small, the correction comes from assuming that the neutron spectrum will take the form

$$\phi(u) \approx \phi_\infty(u) \frac{1}{\Sigma_t(u)(1/\bar{l})} \ , \quad (14)$$

where $\bar{l}$ is the average uncollided path length of a neutron in the slab, also known as the mean chord length. In the thin slab case, $\bar{l} \approx 4V/S$, where $V$ is the volume and $S$ is the surface area of the slab.

The correction is applied to the background cross section by adding an escape cross section, $\bar{\sigma}_{ex}^{(j)} = 1/(N_j \bar{l})$, to it. Thus,

$$\sigma_0^{*,(j)} = \frac{1}{N_j} \sum_{m \neq j} \sigma_t^{(m)} N_m + \frac{1}{N_j \bar{l}} \ , \quad (15)$$

and the weighting spectrum that is used in Eq. (12) becomes $W(u) = \phi_\infty(u)(1/(\sigma_t^{(j)}(u) + \sigma_0^{*,(j)}))$. Note that if $\bar{l}$ is large, the escape cross section effectively becomes zero, and the effect of spatial self-shielding is removed.

## III. PROBLEM CHARACTERIZATION

While standard cross-section processing techniques are sufficient for FW-CADIS to enable the solution of many fixed-source problems with MC methods, this is not always the case. Solutions exhibit very slow MC error reduction in regions that combine space and energy self-shielding. The dose rates predicted by the deterministic code in these regions are a factor of ~10 lower than the MC code.

To illustrate this behavior, throughout this work, we use an idealized geometry with a steel plate embedded in water. This geometry contains the characteristics common to all cases in which the slow error reduction has been observed: a thin slab of material with strong resonances (the iron in steel) inside a thick slab of moderator (water). The resonances in the iron cause self-shielding in energy, and spatial self-shielding comes from having a thin plate of resonance material inside a moderator.

The 53- × 50- × 140-cm stack-up is uniform in the $y$ dimension and in the $x$ dimension, except for the steel plate, which extends from $x = 25$ to 28 cm. The geometry is shown in Fig. 1. The materials vary in the $z$ direction, with units in centimeters, as 0 to 15: source, which is orange; 15 to 30: steel, which is black; 30 to 130: water, which is blue; 130 to 140: air, which is white; and the steel plate extends from 30 to 130 across $x = 25$ to 28. The source is a $^{235}$U fission spectrum that is uniformly distributed throughout a homogenization of water, Zr, and U.

The primary MC code used in this work is MC21 (Ref. 7). Deterministic calculations were executed with PARTISN (Ref. 8), including the three-dimensional forward and adjoint calculations whose results are used to make an importance map for MC21. Unless otherwise noted, PARTISN was used with quadruple range (QR) quadrature sets.[9] The QR sets are defined by the number of polar angles and then the total angles per octant; e.g., QR-8-36 has 8 and 36 angles, respectively. In some cases, MCNP (Ref. 10), which employed manually generated geometric importances, was used for comparison.

MC21 cross-section information comes from ENDF/B-VII. NJOY and TRANSX (Ref. 11) produced the MG cross sections used by PARTISN from ENDF/B-VII data. TRANSX calculates $\sigma_0^{*,(j)}$ using Eq. (15) with a default chord length of 10 000 cm (the infinitely dilute approximation). TRANSX then interpolates the cross sections available in an NJOY library. The user does not have the ability to change $\sigma_0^{*,(j)}$ directly, but they can modify $\bar{l}$. Unless otherwise noted, cross sections were made with a chord length of 10 000 cm. The MCNP calculations used ENDF/B-VII cross-section data as well. All material temperatures were set to 293.6 K.

The results throughout this work are nonthermal (above 0.625 eV) dose rates in millirems per hour unless

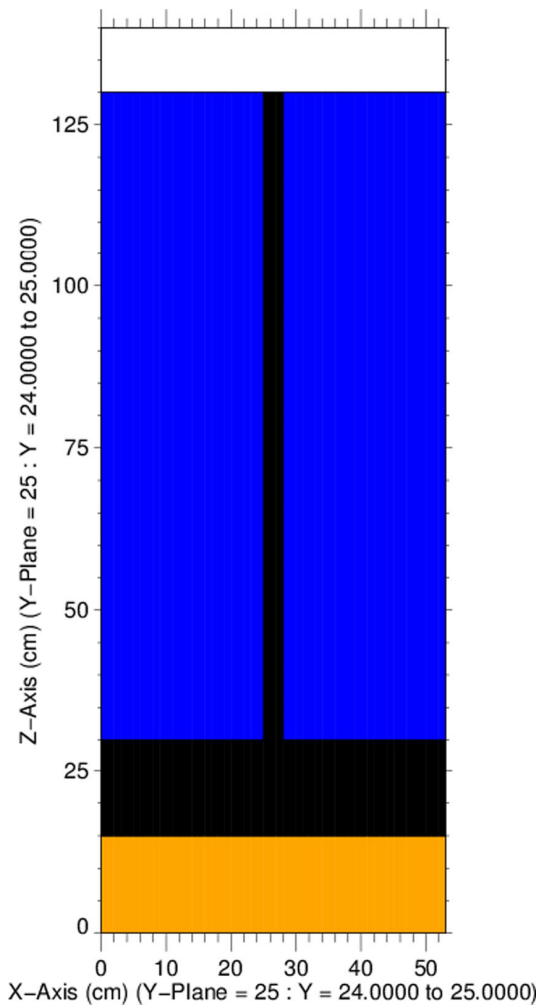

Fig. 1. Geometry (*x-z* slice through *y* = 25 cm).

TABLE I

Base Calculation Parameters

| Variable | PARTISN | MC21 |
|---|---|---|
| Deterministic mesh | 0.5 cm uniform; 0.25 cm in *x* over *x* = 24 to 29 cm | 1 cm uniform |
| Tally mesh | N/A[a] | 1 cm uniform |
| Number of particles | N/A | |
| Energy structure | 58 groups | 27 groups/ continuous |
| Angular quadrature | QR-18-252 | QR-8-36 |
| Scattering order | $P_3$ | $P_3$ |
| Convergence criterion | 0.01 | 0.05 |
| TRANSX settings | Default | Default |
| DCFs | 58 groups | 27 groups |

[a]N/A = not applicable.

otherwise noted. Flux-to-dose conversion factors (DCFs) are applied to the flux via either explicit multiplication (PARTISN) or tally multipliers (MC21) to obtain dose rate. The MCNP results were tallied as nonthermal flux only, so all MCNP comparisons use nonthermal flux rather than dose rate. This work uses 27- and 58-group energy structures, both of which have one thermal group, resulting in 26 and 57 nonthermal groups, respectively.

Table I gives the parameters used in most calculations; differences will be noted as appropriate. The deterministic parameters listed under MC21 were used in PARTISN to create the importance map, and the convergence criterion refers to the PARTISN check on cellwise maximum difference of scalar flux between iterations.[8] Unless otherwise specified, FW-CADIS was set to optimize the calculation of space- and energy-dependent flux everywhere.

The remainder of this section will begin by establishing the fact that MC21, PARTISN, and MCNP agree in the absence of the plate. Next, the manifestation of the problem with the plate is shown. This is followed by a series of investigations attempting to obtain low relative errors (REs) at the exit of the plate (*z* = 130 cm) by refining FW-CADIS parameters, etc. Finally, some calculations using Monaco with Automated Variance Reduction using Importance Calculations[12] (MAVRIC) are presented.

### III.A. Initial Results

To ensure that disagreement in the deterministic and MC results was not coming from simple bulk transmission, a version of the geometry without the plate was run in several codes: PARTISN with a 0.5-cm uniform mesh, QR-8-36, $P_3$, and 27 energy groups; MC21 with FW-CADIS-generated importances using the parameters shown in Table I; and MCNP with cell-based importances set manually using information from the MC21 solution. This will be called the "no-plate" case.

MC21 reports 95% confidence interval (95CI) REs. Figure 2 shows the MC21 nonthermal dose rate and 95CI REs in the absence of a plate as an *x-z* slice through *y* = 25 cm (the center of the geometry). Throughout this work, all MC21 RE color map plots are shown on the same scale to facilitate comparison between plots. In this case's centerline slice, the REs are all under 2% and do not exhibit any odd features.

Comparisons of nonthermal flux between MC21, MCNP using $1 \times 10^9$ particles, and PARTISN are plotted in Figs. 3 and 4. Both plots are along the *z* dimension and use data from the cell covering 26 to 27 cm in the *x* direction and 24 to 25 cm in the *y* direction (through the center of the geometry). Figure 3 contains the nonthermal fluxes from all three codes (refer to the left axis) as well as the ratios of the PARTISN to MC21 codes and the MCNP to MC21 codes (refer to the right axis). All plots containing ratios of nonthermal flux have a reference line of 1.0. Figure 4 shows the REs from the MC codes. The plotted MCNP REs are the MCNP-reported $1\sigma$ RE values multiplied by 1.96 to approximate a 95CI RE.

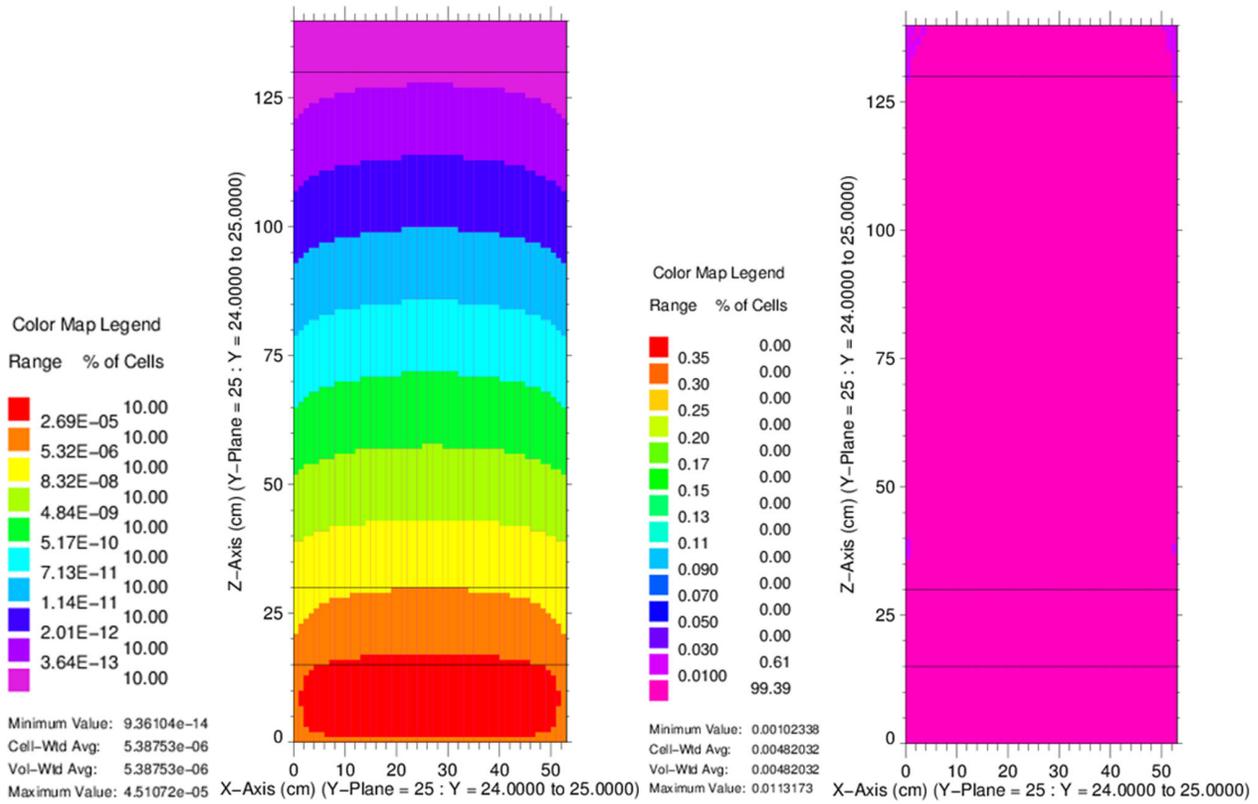

Fig. 2. No-plate MC21 dose rate (left) and 95CI RE (right) (*x-z* slice through *y* = 25 cm).

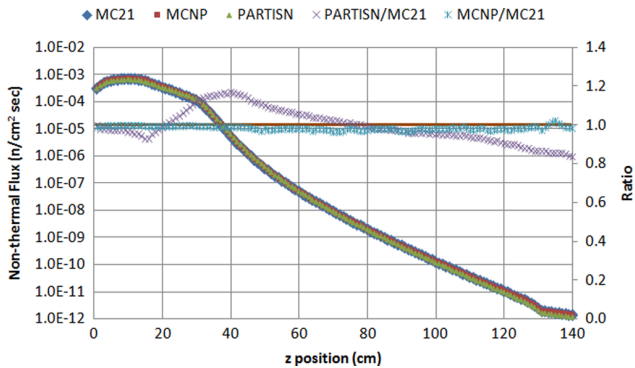

Fig. 3. Plate nonthermal flux (left axis) and method ratios (right axis) down the *x-y* centerline.

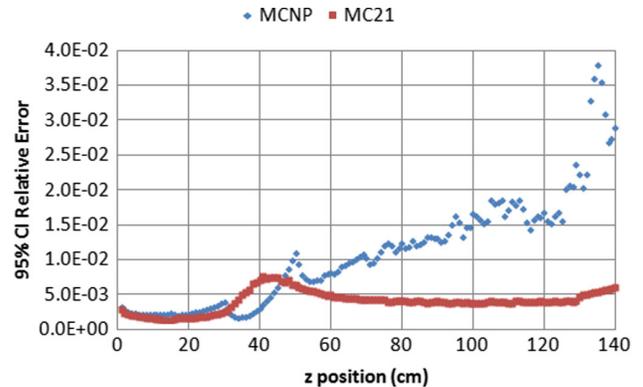

Fig. 4. No-plate 95CI RE for nonthermal flux down the *x-y* centerline.

Both the PARTISN and the MC codes produce very similar results when the steel plate is absent. Compared to the MC results, all of which have <4% 95CI RE, PARTISN's results exhibit a slight bump shortly after the steel-water interface at $z = 30$ cm and subsequently decrease through the remainder of the geometry to a value that is ~16% low. This can be attributed to the approximations employed in the discrete ordinates method. The important point is that through ~4 ft of material, the two methods agree to within ~15%. This is a strong indicator there is no serious flaw in either method's characterization of bulk transport behavior.

Next, calculations with the steel plate were performed. The results, referred to as "base case,"

demonstrate the issues created by the plate. Initially, this problem was solved with the parameters given in Table I using both MC21 and PARTSN. The resulting MC21 nonthermal dose rate and 95CI RE are shown in Fig. 5 as an *x-z* slice through $y = 25$ cm.

The MC21 95CI RE increases substantially as transport continues along the plate, with exit values >25%, but is very low everywhere else in the problem. Increasing the number of particles to $1 \times 10^{11}$ drops the maximum 95CI RE in this slice from 31.3% to 16.4%, but the variance spike through the plate persists. It is clear from the dose rate plot that the highest dose rate on the outside of the shield will follow the plate. For this reason, it is particularly problematic that this is the area with the lowest confidence in the answer.

Figure 6 contains the nonthermal fluxes with the plate from MC21, PARTISN, and MCNP as well as the ratios of PARTISN to MC21 and MCNP to MC21. In the plotted results, MC21 used $1 \times 10^{11}$ particles, while MCNP again used $1 \times 10^{9}$. PARTISN overpredicts MC21 by up to 30% initially, similar to the no-plate case. However, the ratio changes drastically along the length of the plate, such that the exit dose rate is low by a factor of ∼10. This discrepancy is well outside of the MC21 95CI RE bounds, which are shown in Fig. 7 and have a maximum of 5.7% for these data. While MCNP has much higher 95CI RE (also shown in Fig. 7), the general flux trend is close to MC21 throughout most of the geometry, and the MCNP results are therefore taken as confirmation of the MC21 results.

The nonthermal flux in the geometry with a plate does not exhibit good agreement between the MC codes and PARTISN, unlike the no-plate case. These results suggest that the introduction of the plate is the primary cause of PARTISN's significant underprediction of the MC solutions, both in the vicinity of and well outside of the plate. A possible cause of the increase in 95CI RE in the MC21 FW-CADIS solutions in the vicinity of the plate, then, may be using deterministic solutions in the FW-CADIS process that are too inaccurate for this case.

Comparing the spectral tallies from the plate and no-plate cases provides some insight into the source of the deterministic-MC disagreement. In this work, group 1 corresponds to the highest-energy group. The spectra match quite well in the source region ($z = 8.5$ cm) for both geometries, as expected and illustrated in Figs. 8 and 9. Note that the following six plots have lines connecting data points; these are for ease of visualization and are not intended to imply behavior of the data. At $z = 31.5$ cm, the beginning of the large water region, the PARTISN and MC spectra match fairly well for the both the plate and no-plate cases, as can be seen in Figs. 10 and 11. PARTISN is slightly higher in most groups in both cases, but the degree to which it is higher near group 21 for the plate case is larger than the no-plate case.

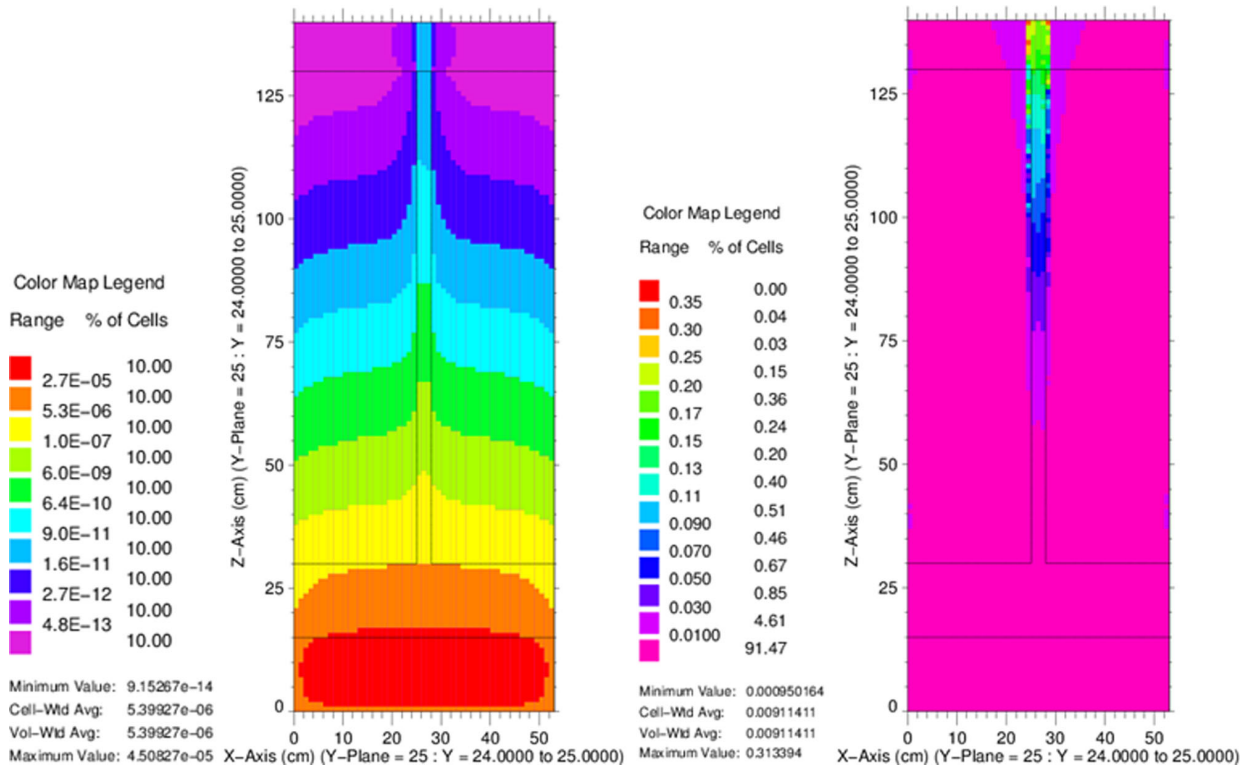

Fig. 5. Plate base FW-CADIS MC21 dose rate (left) and 95CI RE (right) (*x-z* slice through $y = 25$ cm).

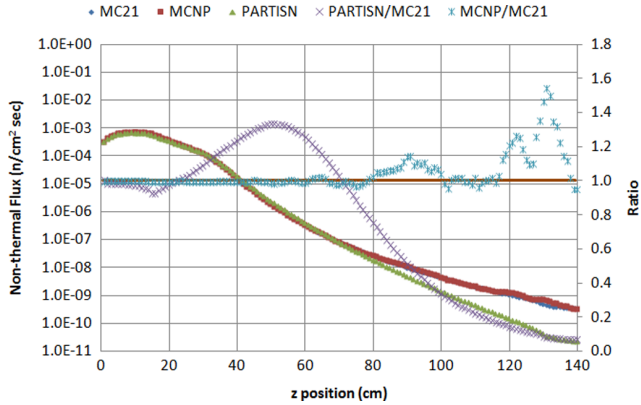

Fig. 6. Plate nonthermal flux (left axis) and method ratios (right axis) down the $x$-$y$ centerline.

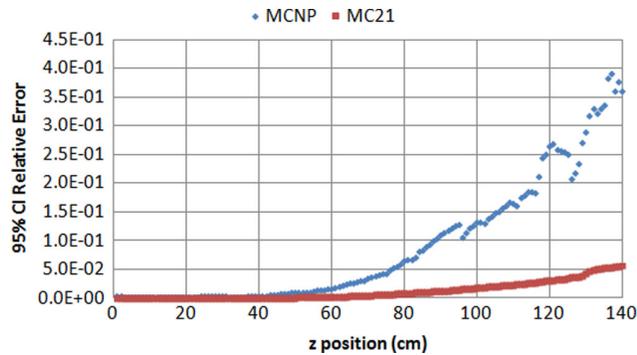

Fig. 7. Plate 95CI RE for nonthermal flux down the $x$-$y$ centerline.

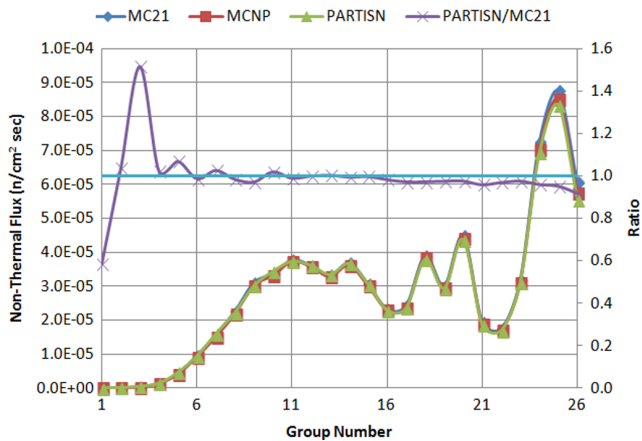

Fig. 8. No-plate flux spectra (left axis) and PARTSN/MC21 (right axis) in source region ($z = 8.5$ cm).

The full extent of the discrepancy, though, appears in the flux spectra as particles exit the plate, $z = 130.5$ cm. The differences between the PARTISN and MC spectra for the no-plate case, shown in Fig. 12, are expected after fission source particles are transported through a thick water/steel stack-up. In the geometry with the plate, however, Fig. 13 shows that MC21's flux is significantly higher than PARTISN's beginning at energies below $\sim 1$ MeV and covering most of the resonances in iron. A very large flux peak in group 17 dominates the nonthermal flux response in both the MC21 and MCNP solutions. The peak exists in the PARTISN solution but underpredicts the MC21 solution substantially. Thus, the MG cross sections fail to capture important physical behavior in the plate, causing PARTISN to drastically underpredict the flux in the resonance energy groups in some areas.

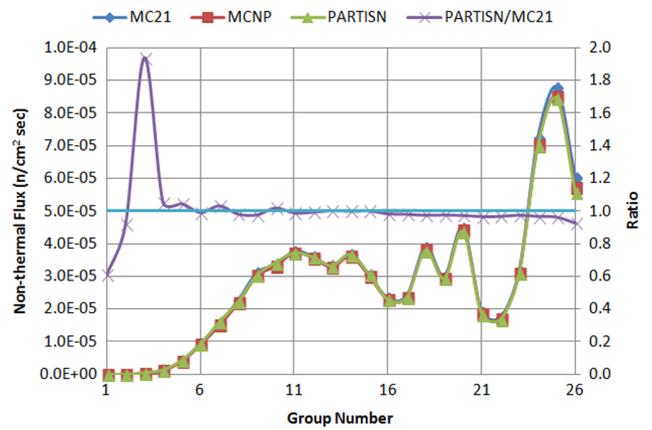

Fig. 9. Plate flux spectra (left axis) and PARTSN/MC21 (right axis) in source region ($z = 8.5$ cm).

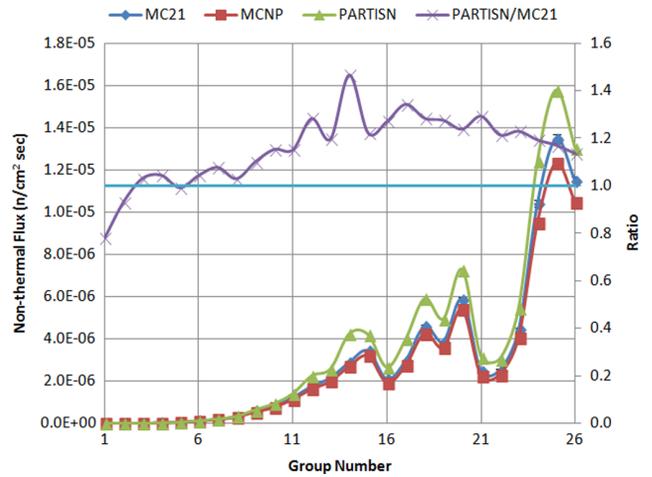

Fig. 10. No-plate flux spectra (left axis) and PARTSN/MC21 (right axis) at start of water region ($z = 31.5$ cm).

We hypothesize that because the degree to which the PARTISN-calculated results are incorrect is large and changes rapidly, the angle-independent importances generated using FW-CADIS are not effective at accelerating MC codes in this region. We think that the ineffective importance map is not caused by a fundamental flaw in FW-CADIS but is a result of using a deterministic solution that is not sufficiently accurate.

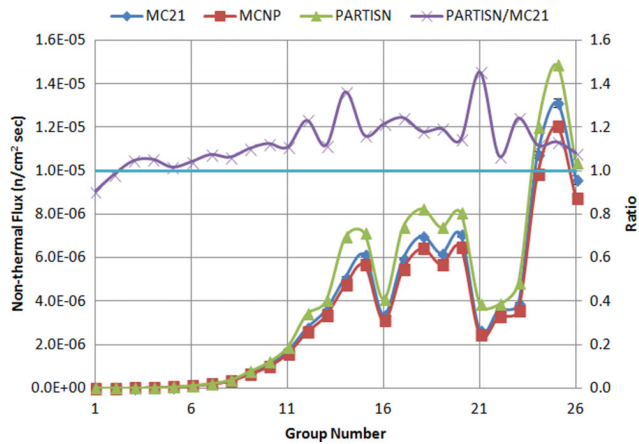

Fig. 11. Plate flux spectra (left axis) and PARTSN/MC21 (right axis) at start of water region ($z = 31.5$ cm).

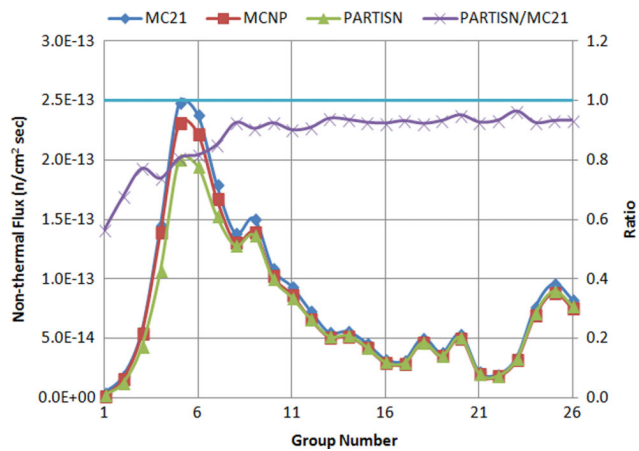

Fig. 12. No-plate flux spectra (left axis) and PARTSN/MC21 (right axis) at start of air region ($z = 130.5$ cm).

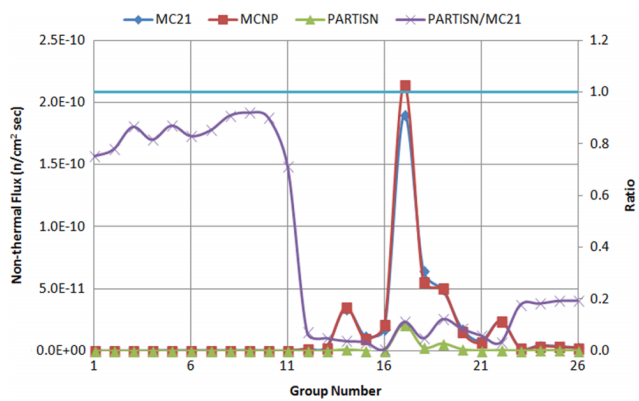

Fig. 13. Plate flux spectra (left axis) and PARTSN/MC21 (right axis) at start of air region ($z = 130.5$ cm).

### III.B. Further Investigation

To check our hypothesis, we tried to improve the deterministic results used in the creation of the importance map. First, the resolution of some of the parameters used to generate importance maps was increased. Next, PARTISN input parameters used to calculate stand-alone PARTISN solutions were varied with the goal of improving the PARTISN solution. The parameters and methods investigated include cross-section generation method, quadrature, scattering expansion, energy structure, and mesh.

For completeness, we also investigated whether CADIS would outperform FW-CADIS in this problem. Finally, other material configurations were used to confirm the pathological case is associated with the physics trends (i.e., energy and spatial self-shielding) rather than this particular configuration.

#### III.B.1. FW-CADIS Options

The energy group structure and mesh resolution were increased to try to change the importance map in a way that will lower the 95CI REs in the plate. An FW-CADIS calculation using 58 energy groups (rather than 27) to create the importance map did not correct the behavior, nor did using a finer importance mesh in the area of the plate. The mesh was refined to be 0.25 cm between $x = 24$ and 29 cm (1 cm on either side of the plate), 0.5 cm in $y$, and 0.5 cm between $z = 30$ and 130 cm (the length of the plate).

Table II summarizes some information about the FW-CADIS cases discussed, where a number appended to a name distinguishes it as using a different number of particles from the standard $1 \times 10^{10}$. Table II includes the number of particles, the CPU-hours, the maximum and average 95CI RE, and the minimum and average figures of merit (FOMs). Here $FOM_{min} = 1/(tR^2_{max})$ and $FOM_{avg} = 1/(tR^2_{avg})$, where $t$ is CPU-hours and $R$ is maximum or average 95CI RE, respectively. The time includes only the MC21 calculation; the importance map generation time was $\sim 1\%$ of the MC21 time and does not significantly impact the FOM. Neither running more particles nor using more energy groups nor using a finer mesh through the plate significantly closed the gap between the maximum and average 95CI REs.

#### III.B.2. Stand-Alone PARTISN

To investigate if the infinitely dilute approximation in the MG cross sections was the source of the inaccuracy of the deterministic solution, a PARTISN calculation was

TABLE II

FW-CADIS Calculation Details

| Case | Particles | CPU-hours | Maximum 95CI RE | Average 95CI RE | Minimum FOM | Average FOM |
|---|---|---|---|---|---|---|
| Base | $1 \times 10^{10}$ | 849.77 | $8.10 \times 10^{-1}$ | $1.02 \times 10^{-2}$ | $1.79 \times 10^{-3}$ | 11.3 |
| Base $1 \times 10^{11}$ | $1 \times 10^{11}$ | 8543.47 | $1.64 \times 10^{-1}$ | $3.28 \times 10^{-3}$ | $4.33 \times 10^{-3}$ | 10.9 |
| 58 groups | $1 \times 10^{10}$ | 954.89 | $5.22 \times 10^{-1}$ | $9.33 \times 10^{-3}$ | $3.85 \times 10^{-3}$ | 12.0 |
| Fine plate | $1 \times 10^{10}$ | 905.34 | 1.55 | $1.71 \times 10^{-2}$ | $4.63 \times 10^{-4}$ | 3.78 |
| Fine $2 \times 10^{11}$ | $2 \times 10^{11}$ | 18 367.67 | $3.43 \times 10^{-1}$ | $4.00 \times 10^{-3}$ | $4.62 \times 10^{-4}$ | 3.40 |

performed with cross sections created using a mean chord length of $4V/S$, the thin slab self-shielding Bondarenko approximation, for all materials. The chord length was calculated separately for the slab of steel extending from $z = 15$ to 30 cm and the thin plate of steel extending through the water and was applied in the appropriate locations. The ratio of dose rate results from $\bar{l} = 10\ 000$ cm (the base PARTISN case) to $\bar{l} = 4V/S$ is shown in Fig. 14.

A significant disparity between the two discrete ordinates solutions is introduced by changing the TRANSX mean chord length from the infinite medium setting to the geometric definition, particularly in and around the steel plate. A change in calculated dose rates of this magnitude is not typical behavior for this parameter change and further suggests that the methods for generating deterministic cross sections may be inadequate for geometries of this configuration. Note that the results obtained using $\bar{l} = 4V/S$ are lower than the base case results. Thus, using a mean chord length of $4V/S$ results in dose rates that diverge even farther from the MC21 results.

Sensitivity of the PARTISN solution to angular quadrature was examined. The ratio of the base case PARTISN dose rate to the dose rate produced using a composite midpoint Gauss (CMG) rule with 591 angles per octant over 40 polar levels is shown in Fig. 15. While differences between the two solutions exist in the vicinity of the plate, they are in good agreement. The maximum ratio in this center slice is 1.29, the minimum is 0.96, and the average is 1.01. PARTISN solution inaccuracy due to inadequate angular quadrature does not appear to explain the much more significant differences observed between PARTISN and MC21.

A case using $P_5$ rather than $P_3$ demonstrated that the Legendre scattering expansion was not the cause of the PARTISN-MC21 mismatch. The maximum ratio through the whole geometry of the $P_3$ case to the $P_5$ case was 1.00, the minimum was 0.964, and the average was 0.994.

Table III summarizes the PARTISN calculation results, where "Chord Length" is the value used in cross-section creation, "Scattering" is the Legendre polynomial order used in the scattering expansion, and the ratios are the minimum, maximum, and average ratios throughout the entire geometry of the base PARTISN case over these cases. The deterministic solution was not changed by using a different mean chord length, a finer angular quadrature, or a higher scattering order in ways that we expect to yield a better importance map.

In our investigations we were not able to generate a PARTISN solution that was accurate in and following the plate of resonance material. Capturing the physics well enough to be able to make a more useful importance map with standard FW-CADIS might be possible if (a) an excellent flux guess is used in the creation of the MG cross sections (i.e., nearly the right answer) and the plate is broken into many small regions each using a unique cross section set or (b) a prohibitive number of energy groups (nearing continuous energy) is used.

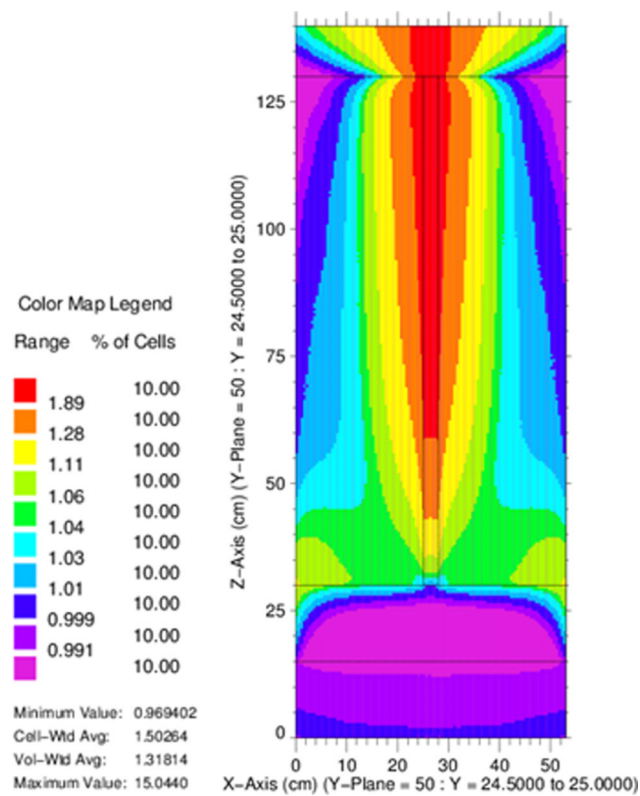

Fig. 14. PARTISN dose rate ratio: $\bar{l} = 10\ 000$ to $\bar{l} = 4V/S$ ($x$-$z$ slice through $y = 25$ cm).

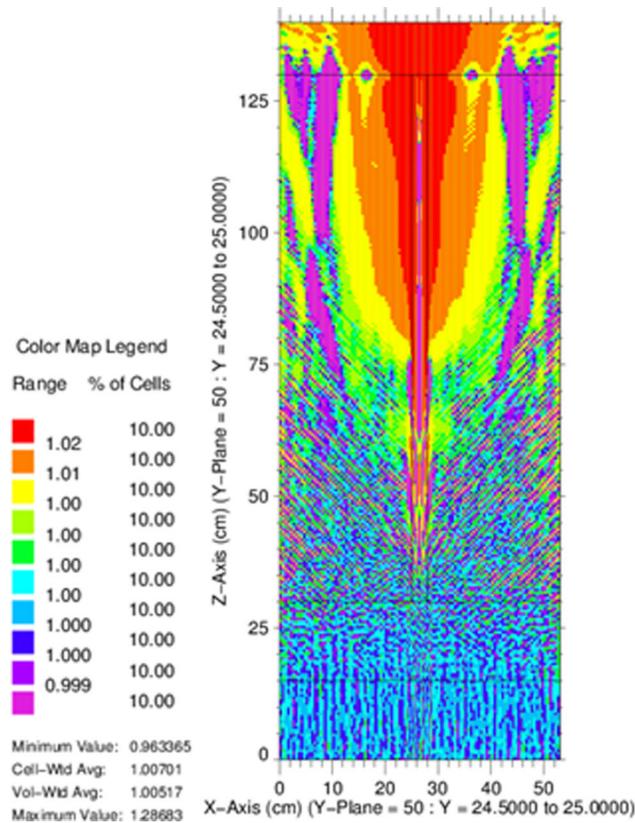

Fig. 15. PARTISN dose rate ratio: QR-18-252 to CMG-591 ($x$-$z$ slice through $y = 25$ cm).

However, we cannot confirm either of these suspicions. The first of these options may not be sufficient because 27 or 58 groups may not capture the physics well enough, even if finely tuned by spatial location. The second option is prohibitively expensive with current computer resources.

### III.B.3. CADIS

CADIS, rather than FW-CADIS, was also used, where the adjoint source location was a 2-cm-thick region following the plate: $x = 25$ to 28 cm, $y = 0$ to 50 cm, and $z = 130$ to 132 cm. Two tests were performed, one in which the adjoint source was set to 1.0 in each energy group, targeting an equally distributed response, and one in which the adjoint source was set to the forward nonthermal flux. The tests produced very similar results, so the test targeting flux will be discussed.

Figure 16 shows the 95CI RE through a centerline slice for the CADIS case. The 95CI RE exiting the plate is now between ∼9% and 20% rather than the base case's 15% and 35%. This is still not an adequate solution for two reasons. First, the 95CI RE exiting the plate remains at least an order of magnitude higher than the surrounding 95CI RE. Second, the quality of the solution in much of the remainder of the problem is poor. For many problems, it is important to know the solution in all locations. Thus, using CADIS instead of FW-CADIS did not resolve the 95CI RE issue.

### III.B.4. Physics Configuration

To ensure that the occurrence of the high 95CI REs is associated with resonance streaming, a case where the plate was chromium, which also has resonances, and one where the plate was air, which does not, was solved. Each of these used the parameters in Table I. Centerline slices of the 95CI REs from the Cr and air cases are shown in Fig. 17. The solutions with Cr have high 95CI RE values. While the 95CI REs are somewhat elevated in the air plate case, they are substantially lower than the cases with resonance materials. The persistently high REs in the air case is a reflection of the difficulty in capturing streaming with an angle-independent importance map. The similarity in 95CI REs between the Cr and Fe cases, in combination with the relative improvement when air is used instead, confirms that the high 95CI REs are associated with the physics involved.

A case with the final air volume changed to polyethylene was also calculated. This served to investigate how adding polyethylene, which will attenuate the neutrons exiting the plate, would impact confidence in the

TABLE III

Summary of PARTISN Stand-Alone Calculations

| Case | Chord Length | Quadrature | Scattering | Maximum Ratio | Minimum Ratio | Average Ratio |
|---|---|---|---|---|---|---|
| 4V/S | 4V/S | QR-18-252 | $P_3$ | 26.23 | 0.93 | 1.34 |
| CMG591 | 10 000 | CMG-591 | $P_3$ | 1.41 | 0.75 | 1.00 |
| P5 | 10 000 | QR-18-525 | $P_5$ | 1.00 | 0.96 | 0.99 |

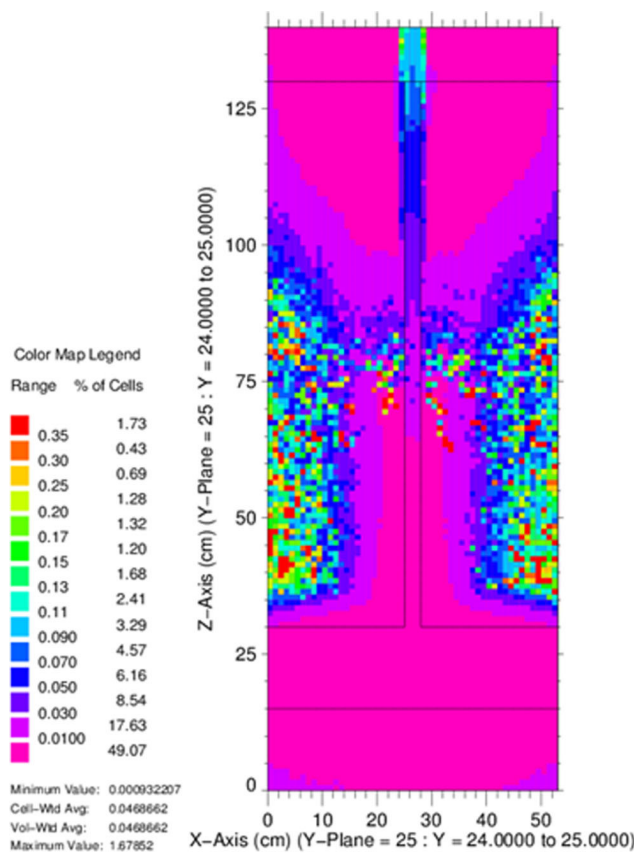

Fig. 16. Flux-targeting CADIS MC21 dose rate 95CI RE (*x-z* slice through *y* = 25 cm).

solution. Figure 18 shows that changing the air to polyethylene mitigates the 95CI RE issue resulting from the standard FW-CADIS calculation. At the interface of the plate and the polyethylene, the 95CI REs are not much lower than when there is only air following. However, as radiation moves into the polyethylene, the 95CI REs decrease significantly compared to the base case. The addition of polyethylene improved the MC21 calculation following the plate.

### III.C. MAVRIC

To investigate the effect of processing the cross sections with more accurate techniques than are available in the version of TRANSX used in this work and to verify that the observed behavior was not due to an error in the implementation of FW-CADIS in the codes used here, we performed independent calculations using MAVRIC. MAVRIC implements FW-CADIS automatically. The forward and adjoint fluxes are calculated by Denovo (Ref. 13), and the MC calculation is carried out by the MG MC code Monaco. All of these codes are part of the SCALE package.[12]

#### III.C.1. Cross Sections in MAVRIC

Within MAVRIC, we elected to use the CENTRM-based cross-section processing path. This path takes infinitely dilute, problem-independent MG data and converts them into problem-dependent MG data by using a weighting function in Eq. (10) that is applicable to the problem being solved. Information about cross-section processing using CENTRM is derived from Ref. 5; consult this reference for more details.

To create the cross sections using CENTRM, first Bondarenko self-shielded cross sections are created from a problem-independent MG library of ENDF/B-VII data. These are then used to create initial self-shielded MG data. The initial self-shielded data, in combination with pointwise (PW) data libraries, are used to compute problem-dependent PW flux spectra. The problem energy range is broken into upper multigroup (UM), PW (the default PW energy range is $E$max = 20 keV to $E$min = $10^{-3}$ eV), and lower multigroup (LM) ranges. Different options can be used to treat the cross sections within each range. Finally, the fluxes for each spatial region are used to compute the final cross sections by numerically integrating Eq. (10) in the PW energy range.

The default CENTRM behavior treats the geometry as an infinite homogeneous medium and performs zonewise homogeneous medium processing in the UM, PW, and LM ranges. The $S_N$ option, which solves the one-dimensional discrete ordinates problem, was used in all three energy ranges. When following the $S_N$ path, MG calculations using standard discrete ordinates techniques are done first in the UM and LM ranges with self-shielded cross sections from the Bondarenko method. The UM and LM calculations are done first to provide sources into the PW range. In the PW energy range, kept as the default in this work, CENTRM solves the $S_N$ equation on a problem-specific PW energy mesh to obtain neutron flux per lethargy as a function of space and direction.

Two calculations were performed with MAVRIC using CENTRM cross-section processing with 200 neutron energy groups (178 nonthermal groups) and one with 27 neutron energy groups (19.45 of which are nonthermal; the additional 0.45 comes from a group that straddles the 0.625-eV cutoff, and the tally function was weighted appropriately). All calculations tallied nonthermal flux, used a $P_3$ scattering expansion, $S_4$ level-symmetric angular quadrature, and a 1-cm uniform mesh with 1/8-cm spacing from $x$ = 24 to 29 cm. In the second 200-group calculation, the plate was broken into nine separate 10-cm segments, each of which uses a different, region-tailored cross-section set to try to obtain more accurate cross sections along the plate. This case is called "200 detailed."

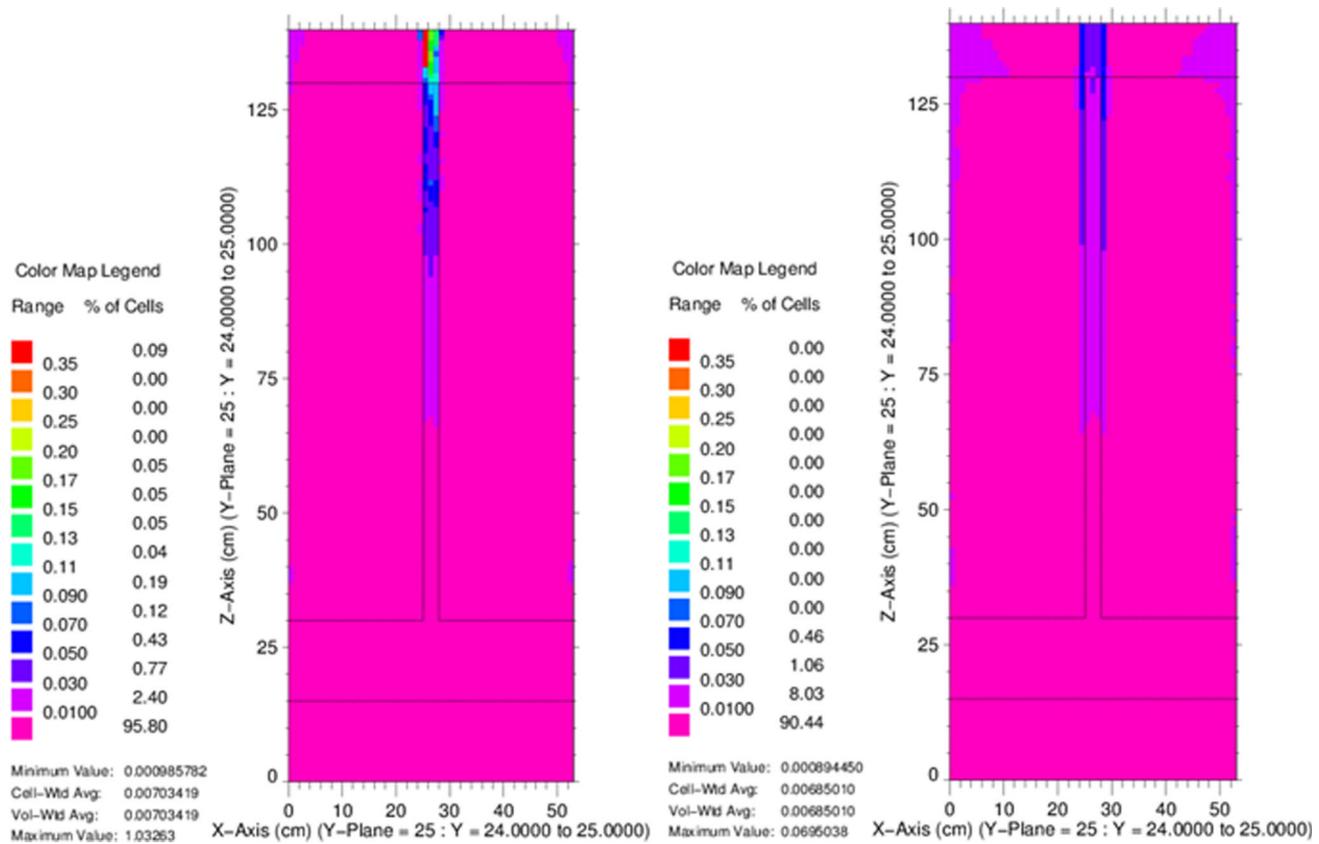

Fig. 17. Chromium plate (left) and air plate (right) FW-CADIS MC21 dose rate 95CI RE ($x$-$z$ slice through $y = 25$ cm).

### III.C.2. MAVRIC Results

Table IV summarizes the MAVRIC results. MAVRIC uses a MG MC code, so the number of groups used in the importance map production is also the number of groups used in tallying the results. The RE results were converted to 95CI RE to be consistent with MC21. Plots of the nonthermal flux and 95CI RE for the 27-group case are in Fig. 19, and those for the 200-group detailed case are in Fig. 20.

The 27-group case misses the streaming behavior and associated high REs entirely. This is likely because (a) the resonance streaming physics is not well captured as a result of the coarse group structure and (b) the directional streaming is not well captured because the importance map is angle independent. The 200-group cases exhibit the streaming and high RE behavior in a manner similar to MC21, confirming that the problem is not caused by our FW-CADIS implementation.

### III.D. Initial Investigation Summary

Section III has discussed a pathological case where high REs are persistent. This behavior is related to space and energy self-shielding caused by a long, thin resonance material embedded in a moderator. Using a finer quadrature, higher scattering order, or a more theoretically accurate mean chord length either did not affect or did not improve PARTISN solutions, which suggests that those changes would not improve importance map generation. Further, using a finer spatial mesh in the problem area or a finer energy group structure to create the importance map did not reduce the high REs. As long as a sufficient number of groups were used, the same patterns were observed when using a more sophisticated cross-section creation procedure. The high RE in and following the plate is present when using different codes and different methods to produce VR parameters and exists when a different resonance material is used.

We suspect that if a sufficiently accurate PARTISN solution were obtained, FW-CADIS would be able to create an importance map that would more substantially enhance MC performance in the plate. However, the computational cost of generating an accurate deterministic solution renders off-the-shelf FW-CADIS inefficient for problems of this nature. This is particularly true since the utility of FW-CADIS comes from its ability to generate highly effective VR parameters with inexpensive and easy-to-acquire deterministic solutions.

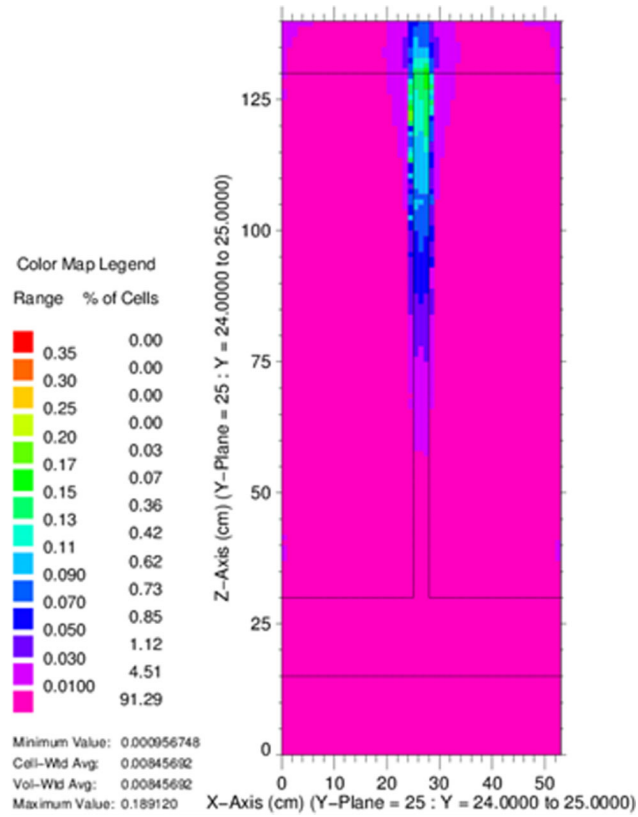

Fig. 18. Polyethylene end block FW-CADIS MC21 dose rate 95CI RE (*x-z* slice through *y* = 25 cm).

TABLE IV

MAVRIC Calculation Results Summary

| Case | Particles | CPU-hours | Maximum 95CI RE | Average 95CI RE | Minimum FOM | Average FOM |
|---|---|---|---|---|---|---|
| 27 group | $1 \times 10^{10}$ | 2086 | $5.56 \times 10^{-3}$ | $1.53 \times 10^{-3}$ | 15.5 | 205 |
| 200 group | $5 \times 10^{9}$ | 1016 | 1.94 | $1.53 \times 10^{-2}$ | $2.62 \times 10^{-4}$ | 4.20 |
| 200 detailed | $9 \times 10^{9}$ | 2098 | 1.93 | $1.22 \times 10^{-2}$ | $1.28 \times 10^{-4}$ | 2.29 |

## IV. METHOD

We have had much success with FW-CADIS, e.g., Ref. 14, but it seems that targeting a uniform response in the way suggested by Wagner et al.[2] does not result in uniform statistical uncertainties for the pathological case discussed in Sec. III. We have determined empirically that applying a correction, which we are calling a Resonance Factor, to the adjoint source specified in FW-CADIS, referred to as $q^{\dagger}_{FWC}$, reduces the RE in the slow error reduction areas:

$$q^{\dagger}(\mathbf{r},E) = \left(\frac{\phi_{res(\sigma_0)}(\mathbf{r},E)}{\phi_{dilute(\sigma_0)}(\mathbf{r},E)}\right)^M q^{\dagger}_{FWC} , \quad (16)$$

where

$M$ = problem-dependent constant typically between 0 and 3

$\phi_{res(\sigma_0)}(\mathbf{r},E)$ = deterministic forward flux calculated with a small background cross section in the resonance material, causing the resonances to have a larger impact

$\phi_{dilute(\sigma_0)}(\mathbf{r},E)$ = deterministic forward flux calculated with a large background cross section in the resonance material, reducing the impact of resonances.

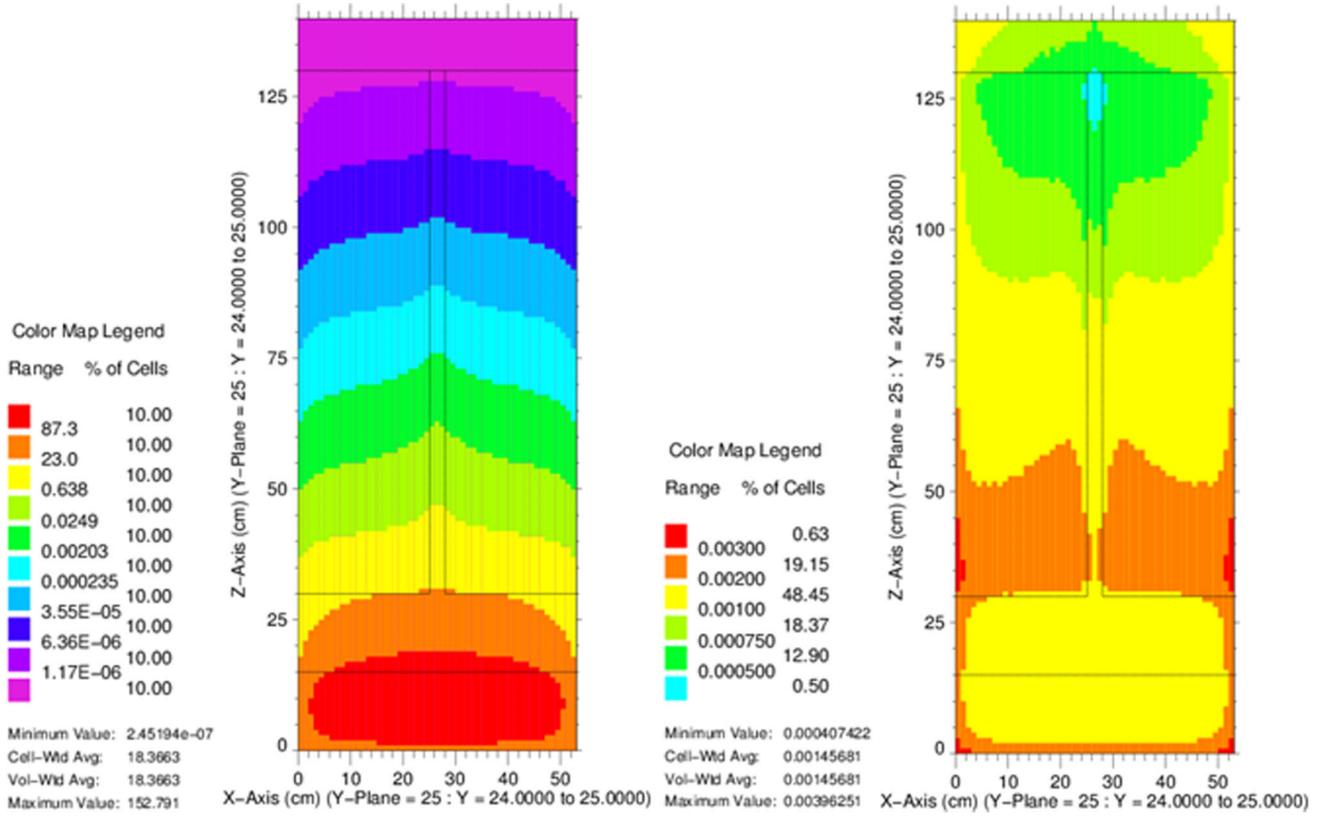

Fig. 19. MAVRIC 27-group nonthermal flux (left) and 95CI RE (right) ($x$-$z$ slice through $y = 25$ cm).

The importance map is created as before, ($imp(\mathbf{r},E) = \phi^\dagger(\mathbf{r},E)/R(\mathbf{r},E)$), but now $\phi^\dagger$ is $\phi^\dagger_{ResFac}$ because it was created with the adjusted adjoint source. Equation (16) results in an adjoint flux that is influenced by the relative importance of the resonance flux compared to the dilute flux. The degree of influence is governed by the magnitude of $M$. Section V studies $M$ for the geometry tested in this work.

To implement the Resonance Factor method, use $\phi_{res(\sigma_0)}(\mathbf{r},E)$ in all locations where $\phi$ would be used in FW-CADIS. Using a small $\sigma_0$ value to generate MG cross sections yields fluxes closer to fluxes created using a best estimate $\sigma_0$ [$\phi_{best(\sigma_0)}(\mathbf{r},E)$, created by tuning the cross section down the length of the plate like what was done in the detailed MAVRIC calculation] than fluxes calculated with a large $\sigma_0$. A more accurate way to implement this new method would be to use $\phi_{best(\sigma_0)}(\mathbf{r},E)$ in all locations where $\phi$ is used in FW-CADIS, but the improvement was not found to be worth the additional cost and complication of computing a third set of cross sections and doing a third forward deterministic calculation.

To illustrate how the correction is applied, consider the example of trying to optimize the space- and energy-dependent flux. In this case, $q^\dagger_{FWC} = 1/\phi$ and results in this corrected response:

$$R = \int_E \int_V \left(\frac{\phi_{res(\sigma_0)}(\mathbf{r},E)}{\phi_{dilute(\sigma_0)}(\mathbf{r},E)}\right)^M \frac{1}{\phi_{res(\sigma_0)}(\mathbf{r},E)} \phi_{res(\sigma_0)}(\mathbf{r},E) dVdE$$

$$\approx \int_E \int_V \left(\frac{\phi_{res(\sigma_0)}(\mathbf{r},E)}{\phi_{dilute(\sigma_0)}(\mathbf{r},E)}\right)^M dVdE . \quad (17)$$

The expanded version of the importance map in the plate then becomes

$$imp(\mathbf{r},E) = \frac{\phi^\dagger_{res(\sigma_0)}(\mathbf{r},E)}{\int_E \int_V \left(\frac{\phi_{res(\sigma_0)}(\mathbf{r},E)}{\phi_{dilute(\sigma_0)}(\mathbf{r},E)}\right)^M dVdE} . \quad (18)$$

## V. RESULTS AND DISCUSSION

The Resonance Factor method as described above was used to solve the geometry containing a plate with varying values of $M$. With the exception of the cross sections used in the calculation of $\phi_{res(\sigma_0)}$ needed in Eq. (16), all other parameter choices are consistent with Table I.

Using the Resonance Factor method significantly improved the behavior of MC21 in the area of the plate. The initial discussion will focus on results using $M = 1$, which gave the best minimum FOM. Figure 21 shows the dose rate 95CI RE from this calculation. The 95CI RE still increases at the end of the plate but largely remains under 10% as opposed to under 20%, as in the base case. Further, the 95CI RE in the rest of the geometry becomes closer to the error in the plate and thus closer to the desired uniform 95CI RE.

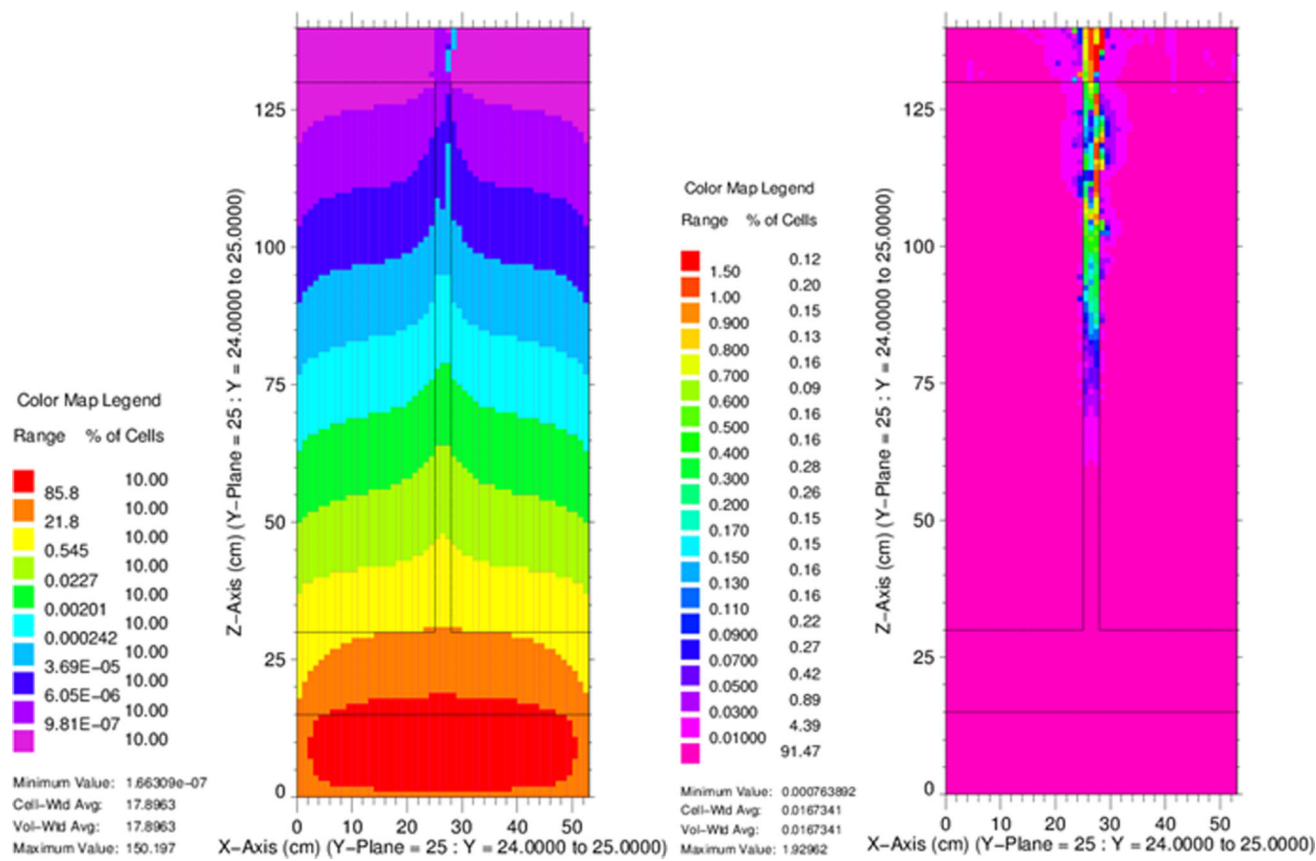

Fig. 20. MAVRIC 200-group, detailed nonthermal flux (left) and 95CI RE (right) ($x$-$z$ slice through $y = 25$ cm).

Table V compares the impact of the new method to the base and best-performing (58-group) FW-CADIS cases across CPU-hours, maximum and average 95CI RE, and minimum and average FOMs. When considering these results, recall that the main goal is to reduce the maximum 95CI RE values that occur following the plate. Having the maximum and average 95CI REs become closer to one another is desirable as well. In this light, it becomes clear that the Resonance Factor is doing a better job of accomplishing this goal. The $FOM_{min}$ is about an order of magnitude higher than the best FW-CADIS case (a factor of 14 higher than the base case), and the $FOM_{min}$ and $FOM_{avg}$ are two orders of magnitude closer together.

The final phase of this study was investigating the impact of the $M$ parameter on performance and 95CI RE. The results as a function of $M$ are shown in Table VI. Within each category, the best value is shown in a bold font and the worst values are underlined.

Figure 22, showing 95CI RE plots for $M = 1 \times 10^{-4}$ and $M = 3.0$, provides a qualitative assessment of the impact of extreme $M$ values. The case with $M = 1 \times 10^{-4}$ provides results that are in close agreement with the base FW-CADIS case (see Fig. 5). This is expected as the maximum normalization is 1.00593 and the minimum is 0.99961, giving an importance map that will not be significantly different from base FW-CADIS. For the $M = 3$ case, the maximum normalization is $5.17641 \times 10^7$, and the minimum is 0.30761. This large range of normalization values provides an inadequate correction to the importance map. The behavior in and just outside of the plate is improved, but at the expense of much of the rest of the problem.

To facilitate comparison across $M$ for the Resonance Factor tests, Table VII shows rank ordering of best to worst of each Resonance Factor case for several of the metrics, focusing on all measures of FOM. The maximum 95CI RE and minimum FOM rankings are nearly the same, as expected. This holds for the minimum 95CI RE/ maximum FOM and average 95CI RE/average FOM as well.

For this problem, the best improvement for the minimum FOM/maximum 95CI RE comes from moderate values of $M$, which correspond to moderately correcting the FW-CADIS adjoint source. Problems with a larger spatial extent and that are more difficult from a MC perspective are likely to benefit from using a larger $M$ value. The more extreme the resonance streaming physics is, the larger the correction needs to be to direct particles to the locations of interest in the problem's geometry.

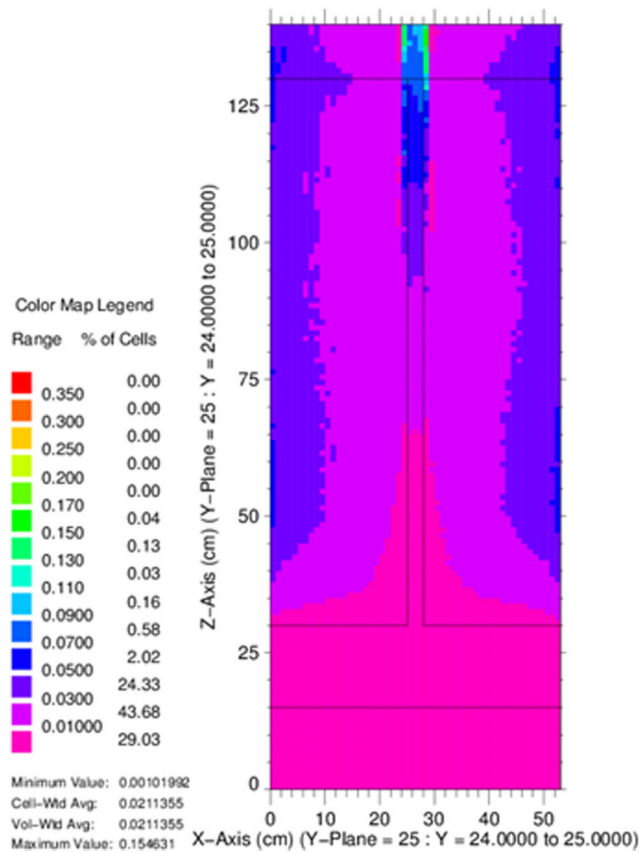

Fig. 21. MC21 dose rate 95CI RE using the Resonance Factor with $M = 1$ (*x-z* slice through $y = 25$ cm).

Another sensible trend is that the average FOM is best for the smallest values of $M$, which give correction values closest to 1.0. Because the FW-CADIS–created importance maps effectively lower the error in nearly all parts of the problem, small values of $M$ will retain this behavior. An unexpected result was that this very small value ($M = 1 \times 10^{-3}$) was useful in that it produced an average FOM value that was better than the best FW-CADIS case (18.24 compared to 11.31, respectively). This is unexpected because the range of normalization values is so small. It seems that through >4 ft of shielding, these small differences compound to have an impact.

Interestingly, the maximum FOM is most improved by the two highest $M$ values. The best maximum FOM is a factor of 3 better than the best FW-CADIS case. This is the case because the Resonance Factor calculations ran more quickly as a result of an improved VR map in difficult parts of the problem, and the minimum 95CI REs are lower as well. These minimum error values occur inside the source region, which is expected, but it is unclear why they are better for the corrected case than the uncorrected.

To summarize the comparisons between the best Resonance Factor and the best FW-CADIS cases, the Resonance Factor was better than FW-CADIS for maximum 95CI RE/minimum FOM and minimum 95CI RE/maximum FOM; FW-CADIS was only better for average 95CI RE/FOM. While we have focused on 95CI RE and FOM, it is worth emphasizing that all Resonance Factor calculations took much less time (~60%) than all FW-CADIS cases. Thus, if trying to reach a certain RE target value for all cells, the Resonance Factor method will accomplish this in less time.

## VI. CONCLUSIONS

We found that the standard FW-CADIS method used with MG cross sections and no angle dependence had difficulty obtaining low MC statistical errors when spatial and energy self-shielding were present in the same area. This issue could not be resolved by refining the angular quadrature, scattering order, spatial mesh, or energy structure; improving the cross-section generation processes; using the Bondarenko resonance correction; or using CADIS alone. To solve this problem, we empirically identified a Resonance Factor, the application of which significantly improves MC error reduction in the challenge area.

These experiences support our theory that the need to use a coarse MG energy structure in deterministic

TABLE V

Comparison of Base FW-CADIS Case, Best FW-CADIS Case, and Best Resonance Factor Case

| Case | CPU-hours | Maximum 95CI RE | Average 95CI RE | Minimum FOM | Average FOM |
|---|---|---|---|---|---|
| Base FW-CADIS | 849.77 | $8.10 \times 10^{-1}$ | $1.02 \times 10^{-2}$ | $1.79 \times 10^{-3}$ | 11.3 |
| 58-group FW-CADIS | 954.89 | $5.22 \times 10^{-1}$ | $9.33 \times 10^{-3}$ | $3.85 \times 10^{-3}$ | 12.0 |
| Resonance Factor $M = 1$ | 534.80 | $2.69 \times 10^{-1}$ | $2.44 \times 10^{-2}$ | $2.58 \times 10^{-2}$ | 3.13 |

TABLE VI

Resonance Factor Results as a Function of $M$

| $M$ | CPU-hours | Maximum 95CI RE | Average 95CI RE | Minimum FOM | Average FOM |
|---|---|---|---|---|---|
| $1.00 \times 10^{-4}$ | 565.80 | $3.56 \times 10^{-1}$ | $9.86 \times 10^{-3}$ | $1.39 \times 10^{-2}$ | 18.18 |
| $1.00 \times 10^{-3}$ | 565.61 | $4.07 \times 10^{-1}$ | $\mathbf{9.85 \times 10^{-3}}$ | $1.07 \times 10^{-2}$ | **18.24** |
| 0.10 | 570.25 | $\mathbf{2.68 \times 10^{-1}}$ | $1.03 \times 10^{-2}$ | $2.44 \times 10^{-2}$ | 16.6 |
| 0.25 | <u>571.09</u> | $3.80 \times 10^{-1}$ | $1.12 \times 10^{-2}$ | $1.21 \times 10^{-2}$ | 13.93 |
| 0.50 | 565.00 | $3.43 \times 10^{-1}$ | $1.39 \times 10^{-2}$ | $1.50 \times 10^{-2}$ | 9.11 |
| 0.75 | 565.00 | $4.02 \times 10^{-1}$ | $1.83 \times 10^{-2}$ | $1.10 \times 10^{-2}$ | 5.31 |
| 1.00 | 534.80 | $2.69 \times 10^{-1}$ | $2.44 \times 10^{-2}$ | $\mathbf{2.58 \times 10^{-2}}$ | 3.13 |
| 1.25 | 522.78 | $2.79 \times 10^{-1}$ | $3.25 \times 10^{-2}$ | $2.46 \times 10^{-2}$ | 1.81 |
| 1.50 | 500.31 | $3.21 \times 10^{-1}$ | $4.23 \times 10^{-2}$ | $1.94 \times 10^{-2}$ | 1.12 |
| 1.75 | 491.78 | $3.97 \times 10^{-1}$ | $5.42 \times 10^{-2}$ | $1.29 \times 10^{-2}$ | 0.69 |
| 2.00 | 483.60 | $7.85 \times 10^{-1}$ | $6.95 \times 10^{-2}$ | $3.36 \times 10^{-3}$ | 0.43 |
| 2.25 | 476.85 | 1.24 | $8.99 \times 10^{-2}$ | $1.37 \times 10^{-3}$ | 0.26 |
| 2.50 | 485.73 | 1.24 | $1.12 \times 10^{-1}$ | $1.33 \times 10^{-3}$ | 0.16 |
| 3.00 | **458.52** | <u>1.81</u> | <u>$1.33 \times 10^{-1}$</u> | $6.62 \times 10^{-4}$ | <u>0.12</u> |

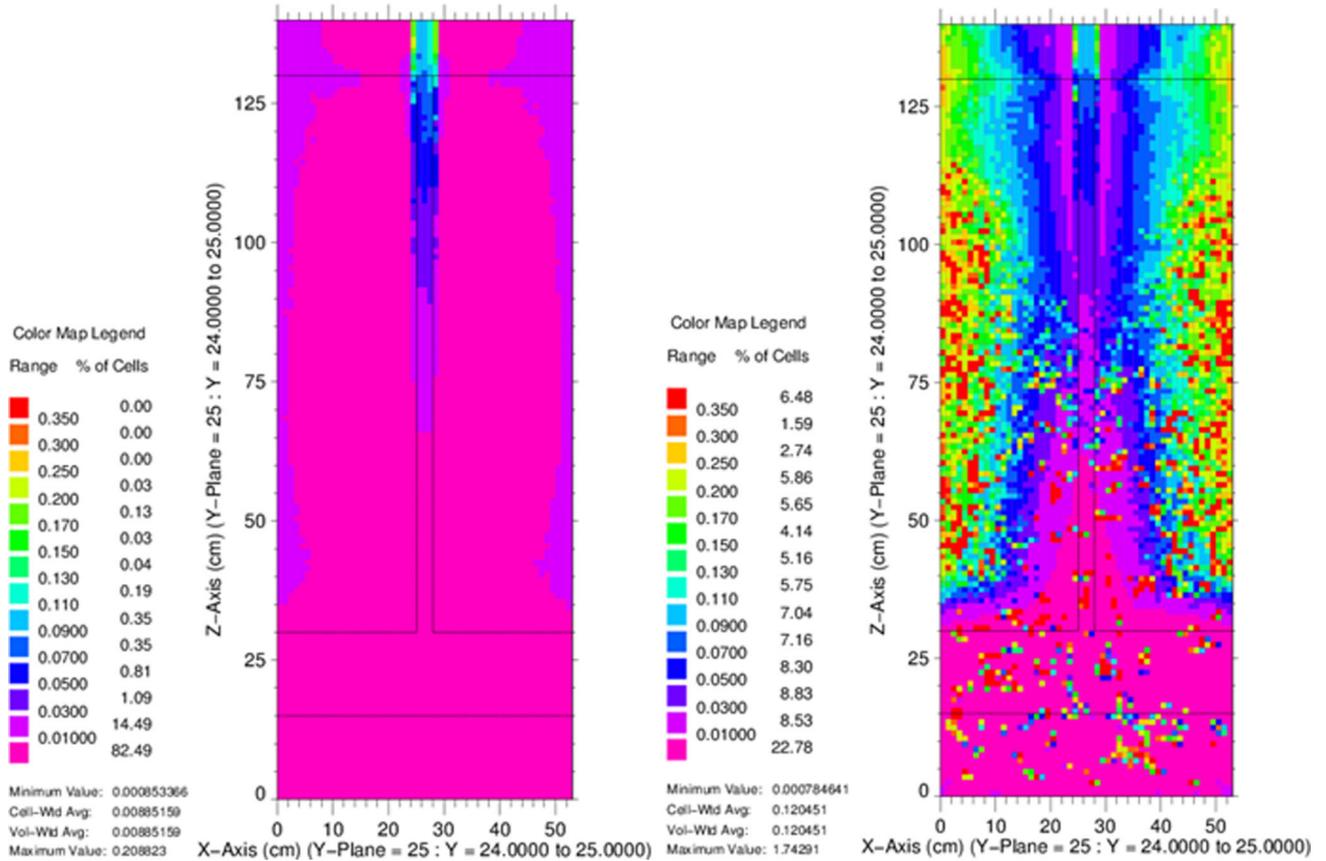

Fig. 22. MC21 dose rate 95CI RE using the Resonance Factor with $M = 1 \times 10^{-4}$ (left) and $M = 3$ (right) (x-z slice through $y = 25$ cm).

TABLE VII

Rank Ordering of Resonance Factor Results from Best to Worst

| Maximum 95CI RE | Minimum FOM | Average FOM | Maximum FOM |
|---|---|---|---|
| $M = 0.1$ | $M = 1.0$ | $M = 1 \times 10^{-3}$ | $M = 3.0$ |
| $M = 1.0$ | $M = 1.25$ | $M = 1 \times 10^{-4}$ | $M = 2.5$ |
| $M = 1.25$ | $M = 0.1$ | $M = 0.1$ | $M = 0.5$ |
| $M = 1.5$ | $M = 1.5$ | $M = 0.25$ | $M = 0.25$ |
| $M = 0.5$ | $M = 0.5$ | $M = 0.5$ | $M = 0.1$ |
| $M = 1 \times 10^{-4}$ | $M = 1 \times 10^{-4}$ | $M = 0.75$ | $M = 1 \times 10^{-3}$ |
| $M = 0.25$ | $M = 1.75$ | $M = 1.0$ | $M = 1 \times 10^{-4}$ |
| $M = 1.75$ | $M = 0.25$ | $M = 1.25$ | $M = 0.75$ |
| $M = 0.75$ | $M = 0.75$ | $M = 1.5$ | $M = 1.0$ |
| $M = 1 \times 10^{-3}$ | $M = 1 \times 10^{-3}$ | $M = 1.75$ | $M = 2.25$ |
| $M = 2.0$ | $M = 2.0$ | $M = 2.0$ | $M = 1.25$ |
| $M = 2.25$ | $M = 2.25$ | $M = 2.25$ | $M = 2.0$ |
| $M = 2.5$ | $M = 2.5$ | $M = 2.5$ | $M = 1.5$ |
| $M = 3.0$ | $M = 3.0$ | $M = 3.0$ | $M = 1.75$ |

calculations combined with the lack of angle dependence in our FW-CADIS implementation causes inadequate VR for resonance energy neutrons that scatter into a flight direction traveling down the plate. We have observed that these particles split many times as they travel down the long, thin plate but do not scatter or otherwise react out of their original energy-angle state. This results in many copies of the same particle being scored on the mesh, which is unnecessary. If the VR were closer to ideal, only unique particles would be scored, and simulating enough unique particles would reduce the variance. The addition of the Resonance Factor ameliorates the copy effect observed with standard FW-CADIS, which implies that the importance map has been improved to better simulate the real physics of the problem.

It is clear that the Resonance Factor method dramatically improves performance in terms of lowering the maximum 95CI RE and raising the minimum FOM. The average performance in the $M = 1.0$ case is not as good as with standard FW-CADIS, but the minimum and average FOMs become closer together. The Resonance Factor is only needed in cases when the answer in and/or following the self-shielded location must be accurate; otherwise, simply use FW-CADIS. Thus, the new method significantly improves the ability to solve these pathological cases.

There are several important things to note about the Resonance Factor method compared to FW-CADIS. It requires three deterministic calculations instead of two: first, a forward calculation using the small $\sigma_0$; second, a forward calculation using the large $\sigma_0$; and third, an adjoint calculation using the small $\sigma_0$. Also, two sets of MG cross sections are needed instead of just one. Finally, an understanding of how to produce these cross sections is required.

This process requires more work than the standard FW-CADIS method. However, without these corrections the amount of computer time needed to obtain sufficiently low statistical errors in and following the resonance material could be prohibitive. We recommend that this method be used when the accuracy of the solution is important in the presence of combined space and energy self-shielding.

REFERENCES

1. D. E. PEPLOW, "Comparison of Hybrid Methods for Global Variance Reduction in Shielding Calculations," *Trans. Am. Nucl. Soc.*, **107**, 512 (2012).

2. J. C. WAGNER, E. D. BLAKEMAN, and D. E. PEPLOW, "Forward-Weighted CADIS Method for Global Variance Reduction," *Trans. Am. Nucl. Soc.*, **97**, 630 (2007).

3. I. I. BONDARENKO, "Group Constants for Nuclear Reactor Calculations," Consultants Bureau Enterprises (1964).

4. M. A. COOPER and E. W. LARSEN, "Automated Weight Windows for Global Monte Carlo Particle Transport Calculations," *Nucl. Sci. Eng.*, **137**, 1 (2001); http://dx.doi.org/10.13182/NSE00-34.

5. M. L. WILLIAMS, "Resonance Self-Shielding Methodologies in SCALE 6," *Nucl. Technol.*, **174**, 149 (2011); http://dx.doi.org/10.13182/NT09-104.

6. "NJOY99.0: Code System for Producing Pointwise and Multigroup Neutron and Photon Cross Sections from ENDF/B Data," Technical Report RSICC Code Package PSR-480, Los Alamos National Laboratory (2000).

7. T. M. SUTTON et al., "The MC21 Monte Carlo Transport Code," *Proc. Mathematics and Computations and Supercomputing in Nuclear Applications (M&C+SNA 2007)*, Monterey, California, April 15–19, 2007, American Nuclear Society (2007).

8. R. E. ALCOUFFE et al., "PARTISN: A Time-Dependent, Parallel Neutral Particle Transport Code System," LA-UR-08-07258, Los Alamos National Laboratory (2009).

9. I. K. ABU-SHUMAYS and C. E. YEHNERT, "Comparison of Hybrid Methods for Global Variance Reduction in Shielding Calculations," *Proc. Topl. Mtg. Radiation Protection and Shielding*, Falmouth, Massachusetts, April 21–25, 1996, American Nuclear Society (1996).

10. X-5 MONTE CARLO TEAM, "MCNP—A General Monte Carlo N-Particle Transport Code, Version 5, Volume 1: Overview and Theory," LA-UR-03-1987, Los Alamos National Laboratory (2008).

11. R. E. MacFARLANE, "TRANSX 2: A Code for Interfacing MATXS Cross-Section Libraries to Nuclear Transport Codes," LA-12312-MS, Los Alamos National Laboratory (1992).

12. "SCALE Version 6.1, Section S6," ORNL/TM-2005/39, Radiation Safety Information Computational Center, Oak Ridge National Laboratory (2011).

13. T. M. EVANS et al., "Denovo: A New Three-Dimensional Parallel Discrete Ordinates Code in SCALE," *Nucl. Technol.*, **171**, 171 (2010); http://dx.doi.org/10.13182/NT10-5.

14. R. N. SLAYBAUGH and S. C. WILSON, "Deterministic Parameter Study for Fixed-Source Calculations Using FW-CADIS," *Trans. Am. Nucl. Soc.*, **108**, 441 (2013).